\documentclass{article}
\usepackage{amsmath}
\usepackage{amsfonts}
\usepackage{amsthm}
\usepackage{enumerate}

\newtheorem{thm}{Theorem}[section]
\newtheorem{prop}[thm]{Proposition}
\newtheorem{cor}[thm]{Corollary}

\newtheorem{lem}[thm]{Lemma}

\theoremstyle{remark}
\newtheorem{rmk}[thm]{Remark}

\begin{document}

\title{The Origin of the Beauregard--Suryanarayan Product on Pythagorean Triples}

\author{Shaul Zemel}

\maketitle

\section*{Introduction}

Pythagorean triples are one of the oldest and most studied objects in elementary mathematics. They appear as the basic examples in many mathematical questions, and trying to generalize them has led to the development of several fascinating fields in mathematics over the centuries---one such direction was the research revolving around, and stemming from, Fermat's Last Theorem. Several papers containing also historical perspectives already appear in the literature, including some remarks in the rather recent reference \cite{[M]}. Since the literature about this subject is vast, any attempt to treat a reasonable amount of references will do injustice with many unmentioned authors. We therefore consider here only those papers that are closely related to our point of view, leaving the more comprehensive lists of references to more historical surveys.

Normalizing a Pythagorean triple according to the largest number of the triple produces elements of the complex unit circle, thus giving one type of product on Pythagorean triples. Explicitly, one can identify the primitive Pythagorean triple $(a,b,c)$ with the element $\frac{a}{c}+\frac{b}{c}\sqrt{-1}$ of the intersection of $\mathbb{Q}(\sqrt{-1})$ with the complex unit circle, and the resulting multiplication on Pythagorean triples (or their classes modulo scalar multiplications), considered by, e.g., \cite{[E]} and \cite{[T]}, is given by \[(a,b,c)\cdot_{i}(f,g,h):=(af-bg,ag+bf,ch)/\gcd\{af-bg,ag+bf\}.\] There are several different normalizations for this product, depending on the signs of the entries $(a,b,c)$ and the relations to primitivity (i.e., \cite{[E]} considers only triples with positive entries, and normalizes the product accordingly, essentially by working modulo the order 4 subgroup of the unit circle that is generated by $\sqrt{-1}$). On the other hand, the reference \cite{[BS]} defines another product with respect to which the set of Pythagorean triples is a commutative semi-group, one that is inherently based on one of the smaller coordinates of the triple. In this paper we shall use the normalization that does not preserve primitivity and takes into account the possible signs, which is given by
\begin{equation}
(a,b,c)\cdot(f,g,h):=(af,bh+cg,bg+ch). \label{BSprod}
\end{equation}
Several results are proved about the structure of this product in that reference, with additional inspection appearing in \cite{[M]}. Both references also treat the basic involution arising from interchanging the two coordinates. It is this product that lies in the heart of our presentation.

\smallskip

In this note we explain the geometric origin of the product defined in \cite{[BS]} as in Equation \eqref{BSprod}, also giving more conceptual reasonings for the results of that reference. The idea is that apart from quadratic fields, there is another 2-dimensional extension $\mathbb{E}$ of the rationals that involves no nilpotent elements (i.e., a \emph{quadratic \'{e}tale algebra} over $\mathbb{Q}$). As quadratic fields appear in the theory of anisotropic rational binary quadratic forms, this algebra $\mathbb{E}$ is the corresponding counterpart for the isotropic quadratic form, also known as the \emph{hyperbolic plane}. The resulting investigation of $\mathbb{E}$ yields naturally the product from \cite{[BS]} and Equation \eqref{BSprod} on Pythagorean triples, and produces directly the structure theorems from that reference. Using the natural involution on Pythagorean triples (which is a Cayley transform in this point of view), we obtain new insights into the relation between right triangles with integral edges that are almost isosceles and Pell's equation with discriminant 2. More precisely, we determine all the integral right triangles with leg difference 1 in terms of coefficients of consecutive solutions of this Pell's equation (or equivalently of consecutive powers of the fundamental unit $1+\sqrt{2}$ of the quadratic field $\mathbb{Q}(\sqrt{2})$). We extend this result to determining Pythagorean triples in which more general linear combinations of the entries are assumed to be small. More explicitly, we show that for every $q\geq2$, the Pythagorean triples $(a,b,c)$ in which the difference $b-qa$ (resp. $c-qa$) is $\pm1$ are all given in terms of consecutive solutions of Pell's equation with the discriminant $q^{2}+1$ (resp. $q^{2}-1$), namely powers of $q+\sqrt{q^{2}+1}$ (resp. $q+\sqrt{q^{2}-1}$), with an exception (which we also determine) for $c-2a$ arising from integral elements of norm $-2$ in $\mathbb{Q}(\sqrt{3})$. For differences like $a-qb$ or $c-qb$ one can apply the natural involution, and obtain similar results.

\smallskip

This paper is divided into 6 sections. Section \ref{QuadForms} presents the rational binary quadratic forms. Section \ref{PythTrip} shows how the algebra $\mathbb{E}$ arises in this context and relates it to the Pythagorean triples, while more explicit expressions and the structure theorem are given in Section \ref{Explicit}. Section \ref{Inv} introduces the natural involution on Pythagorean triples. Section \ref{QuadUnits} relates differences in Pythagorean triples (as ``almost fixed points'' of the involution) to Pell's equation in $\mathbb{Q}(\sqrt{2})$. Finally, Section \ref{LinComb} investigates the differences $b-qa$ and $c-qa$ from Pythagorean triples in terms of units in orders in the fields $\mathbb{Q}(\sqrt{q^{2}\pm1})$, and explains the difficulties that one encounters in the consideration of general linear relations like $pb-qa$ and $pc-qa$.

\smallskip

Many thanks are due to the anonymous referee, for the suggestions how to improve the presentation in this paper.

\section{Quadratic Fields and Binary Quadratic Forms \label{QuadForms}}

We begin by presenting the setup in which we shall later view the group structure on Pythagorean triples. Consider a non-degenerate quadratic space $V$ of dimension 2 over $\mathbb{Q}$, and recall that the determinants of all the possible Gram matrices for $V$ differ from one another by invertible rational squares, i.e., by elements of the group $(\mathbb{Q}^{\times})^{2}$. Hence $V$ determines a \emph{determinant} in the quotient group $\mathbb{Q}^{\times}/(\mathbb{Q}^{\times})^{2}$, and since the dimension is 2 we define the \emph{discriminant} $d$ of $V$ to be \emph{minus} that determinant.

We can take an orthogonal basis for $V$, rendering it isomorphic to the space $\mathbb{Q}^{2}$ with the quadratic form $q(x,y)=ax^{2}-ady^{2}$ for some pre-image of $d$ in $\mathbb{Q}^{\times}$ (which we also denote by $d$). The orthogonal group $\operatorname{O}(V)$, the special orthogonal group $\operatorname{SO}(V)$, and the kernel $\operatorname{SO}^{1}(V)$ of the spinor norm from $\operatorname{SO}(V)$ are all unaffected by rescaling of $V$ (and so is the discriminant in even dimensions), so that we may consider just the quadratic form $q(x,y)=x^{2}-dy^{2}$ for investigating them.

Assume first that $d$ is not a rational square, and let $\mathbb{K}$ be the (real or imaginary) quadratic field $\mathbb{Q}(\sqrt{d})$, whose unique non-trivial Galois automorphism over $\mathbb{Q}$ we denote by $\sigma$. We can identify an element $(x,y)$ of $\mathbb{Q}^{2}$ with the element $z=x+y\sqrt{d}$ of $\mathbb{K}$, and then $q$ is just the norm map $N^{\mathbb{K}}_{\mathbb{Q}}$ from ${\mathbb{K}}$ to ${\mathbb{Q}}$. There is an operation of the group $\mathbb{K}^{\times}$ on the rational vector space $\mathbb{K}$ that preserves the quadratic form (i.e., the norm), in which $\xi\in\mathbb{K}^{\times}$ multiplies $z\in\mathbb{K}$ by $\xi$ but divides it by its Galois conjugate $\xi^{\sigma}$. We therefore obtain a map from $\mathbb{K}^{\times}$ into the orthogonal group $\operatorname{O}(V)=\operatorname{O}(\mathbb{K})$. A classical result from the theory of quadratic forms, which is easily verified directly, is the following.
\begin{prop}
The image of $\mathbb{K}^{\times}$ in $\operatorname{O}(\mathbb{K})$ is precisely the special orthogonal group $\operatorname{SO}(\mathbb{K})$, and the kernel of that map is $\mathbb{Q}^{\times}$. The spinor norm of the image of an element $\xi\in\mathbb{K}^{\times}$ inside $\operatorname{SO}(\mathbb{K})$ is the image of $N^{\mathbb{K}}_{\mathbb{Q}}(\xi)$ in $\mathbb{Q}^{\times}/(\mathbb{Q}^{\times})^{2}$. \label{GSpin}
\end{prop}
Proposition \ref{GSpin} can be described as a short exact sequence \[1\to\mathbb{Q}^{\times}\to\mathbb{K}^{\times}\to\operatorname{SO}(\mathbb{K})\to1\] of Abelian groups, or equivalently by the statement that $\mathbb{K}^{\times}$ is the \emph{GSpin group} of the quadratic space $\mathbb{K}$ (hence of $V$), which in symbols is written as $\operatorname{GSpin}(V)=\mathbb{K}^{\times}$. The full group $\operatorname{O}(V)$ is obtained by adding the action of $\sigma$, so the appropriate cover $\operatorname{Gpin}(V)$ of $\operatorname{O}(V)$ is the semi-direct product in which the Galois group $\operatorname{Gal}(\mathbb{K}/\mathbb{Q})$, of order 2, operates on $\mathbb{K}^{\times}$.

Denote by $\mathbb{K}^{1}$ the subgroup of $\mathbb{K}^{\times}$ consisting of those elements $\xi\in\mathbb{K}^{\times}$ with $N^{\mathbb{K}}_{\mathbb{Q}}(\xi)=1$. Combining the second assertion of Proposition \ref{GSpin} with the obvious fact that $N^{\mathbb{K}}_{\mathbb{Q}}(r)=r^{2}$ for any $r\in\mathbb{Q}^{\times}$ produces the following classical result.
\begin{cor}
The group $\mathbb{K}^{1}$ maps onto the spinor kernel $\operatorname{SO}^{1}(\mathbb{K})$, with kernel $\{\pm1\}$. \label{Spin}
\end{cor}
The short exact sequence associated with Corollary \ref{Spin} is
\begin{equation}
1\to\{\pm1\}\to\mathbb{K}^{1}\to\operatorname{SO}^{1}(\mathbb{K})\to1, \label{shexseq}
\end{equation}
so that $\mathbb{K}^{1}$ is the \emph{spin group} of $V=\mathbb{K}$, i.e., $\operatorname{Spin}(V)=\mathbb{K}^{1}$. With our normalization of $V=\mathbb{K}$ with $q=N^{\mathbb{K}}_{\mathbb{Q}}$, an element of order 2 and spinor norm 1 in $\operatorname{Gpin}(V)$ is the combination of $\sqrt{d}\in\mathbb{K}^{\times}$ with the action of $\sigma$ (since $\sqrt{d}^{\sigma}=-\sqrt{d}$, this action sends $z=x+y\sqrt{d}$ to $-z^{\sigma}=-x+y\sqrt{d}$, indeed representing the reflection in the vector $1\in\mathbb{K}$ with $q(1)=1$). The pin group $\operatorname{pin}(\mathbb{K})$ is therefore a semi-direct product similar to $\operatorname{Gpin}(\mathbb{K})$. However, as this group does vary with different rescalings of $V$, we shall not consider it further.

\begin{rmk}
The group $\mathbb{K}^{1}$ from Corollary \ref{Spin} is \emph{not} the group of units in the ring of integers of $\mathbb{K}$, but is much larger. Indeed, in case $\mathbb{K}$ has class number 1, every prime number $p$ that splits in $\mathbb{K}$ produces an element $\frac{\pi}{\pi^{\sigma}}$ of $\mathbb{K}^{1}$ for some generator $\pi$ of a prime ideal in $\mathbb{K}$ lying over $p$, and it is clear (by considering ideals) that all these generators are independent in the group $\mathbb{K}^{1}$ (also modulo units in the ring of integers). Replacing $\pi$ by another generator simply multiplies the generator by a unit, and taking the generator for the other ideal lying over $p$ replaces the corresponding element of $\mathbb{K}^{1}$ by its inverse (again, perhaps up to units). \label{strucK1}
\end{rmk}

\section{Pythagorean Triples \label{PythTrip}}

Let $(a,b,c)$ be a \emph{Pythagorean triple}, i.e., a triplet of integers satisfying the Pythagorean equality $a^{2}+b^{2}=c^{2}$, and assume that $a\neq0$. Then the rational numbers $\alpha=\frac{c}{a}$ and $\beta=\frac{b}{a}$ satisfy the equality $\alpha^{2}-\beta^{2}=1$, which is the norm 1 condition associated with the case of the trivial discriminant $d=1$. While $\mathbb{Q}(\sqrt{1})$ is no longer a quadratic field, we recall that for a non-square $d$ the quadratic field $\mathbb{K}$ can be presented as the quotient $\mathbb{Q}[X]/(X^{2}-d)$ of the polynomial ring $\mathbb{Q}[X]$. This construction produces, with $d=1$, a ring called the \emph{split quadratic \'{e}tale algebra} $\mathbb{E}=\mathbb{Q}[\varepsilon]$, with $\varepsilon$ being a formal element satisfying $\varepsilon^{2}=1$. The basic properties of $\mathbb{E}$ that we shall need are given in the following lemma, whose proof is simple and straightforward.
\begin{lem}
The identity $1\in\mathbb{E}$ is the sum of the complementary idempotents $\frac{1+\varepsilon}{2}$ and $\frac{1-\varepsilon}{2}$. The map sending a pair $(r,s)$ of rational numbers to the element $r\frac{1+\varepsilon}{2}+s\frac{1-\varepsilon}{2}$ is an isomorphism of rings from $\mathbb{Q}\times\mathbb{Q}$ onto $\mathbb{E}$. The norm $N^{\mathbb{E}}_{\mathbb{Q}}$ takes the element $\frac{r+s}{2}+\frac{r-s}{2}\varepsilon\in\mathbb{E}$ that is associated with the pair $(r,s)$ to the product $rs$. \label{Eprop}
\end{lem}
The Galois automorphism $\sigma$ of a quadratic field $\mathbb{K}$ has an analogue for the algebra $\mathbb{E}$. Indeed, as its action takes, in the field case, the image $\sqrt{d}$ of the element $X\in\mathbb{Q}[X]$ to its additive inverse, here it is defined by $\varepsilon\mapsto-\varepsilon$. It therefore operates by interchanging the two factors of the isomorph $\mathbb{Q}\times\mathbb{Q}$ of $\mathbb{E}$ from Lemma \ref{Eprop}, and we denote it by $\iota$. The fact that $\mathbb{E}$ is not a field but rather contains non-zero non-invertible elements of norm 0 means that as a quadratic space it contains \emph{isotropic vectors}---indeed, it is the quadratic space classically known as the \emph{(rational) hyperbolic plane}.

Now, Proposition \ref{GSpin} and Corollary \ref{Spin} do not depend on $\mathbb{K}$ being a field or the non-triviality of $d$, so that they extend to the case $d=1$ with $\mathbb{K}$ replaced by $\mathbb{E}$. This also holds for the short exact sequences and the statements about the groups $\operatorname{Gpin}(V)$ and $\operatorname{pin}(\mathbb{E})$, with $\iota$ instead of $\sigma$. Writing this case more explicitly, the associated quadratic space (i.e., the hyperbolic plane) is generated by two isotropic vectors $\zeta$ and $\omega$ pairing to 1 with one another, and any element of $\operatorname{SO}(V)$ multiplies $\zeta$ by a non-zero rational number $c$ while dividing $\omega$ by $c$. Such an element has spinor norm $c$ in $\mathbb{Q}^{\times}/(\mathbb{Q}^{\times})^{2}$, and for obtaining $\operatorname{O}(V)$ we add the element inverting $\zeta$ and $\omega$. By identifying $V$ with $\mathbb{E}$, we can take $\zeta$ and $\omega$ to be the idempotents $\frac{1\pm\varepsilon}{2}$.

We can now prove the first result from \cite{[BS]} and \cite{[M]} about the product from Equation \eqref{BSprod}. Observe that the homogeneity of this product implies that it descends (like the product from \cite{[E]}) to a well-defined operation on the equivalence classes of Pythagorean triples under the relation identifying a Pythagorean triple with all its scalar multiples.
\begin{prop}
Consider the set of classes of Pythagorean triples $(a,b,c)$ with $a\neq0$ modulo scalar multiplication, with the product resulting from Equation \eqref{BSprod}. The resulting structure is of a group that is isomorphic to $\mathbb{Q}^{\times}$. \label{projgp}
\end{prop}

\begin{proof}
We have already seen that to a Pythagorean triple $(a,b,c)$ with $a\neq0$ we can associate the element $\alpha+\beta\varepsilon\in\mathbb{E}^{1}$ with $\alpha=\frac{c}{a}$ and $\beta=\frac{b}{a}$, and it is clear that scalar multiples of the same Pythagorean triple correspond to the same element of $\mathbb{E}^{1}$. Conversely, any element $\alpha+\beta\varepsilon\in\mathbb{E}^{1}$ is obtained in this way from an integral Pythagorean triple, as is easily seen by taking $a$ such that $a\alpha$ and $a\beta$ are integral and setting $b=a\beta$ and $c=a\alpha$. Moreover, if another Pythagorean triple, say $(f,g,h)$ with $f\neq0$, yields the element $\gamma+\delta\varepsilon$ of $\mathbb{E}^{1}$, then it is easily verified that their product $(af,bh+cg,bg+ch)$ from Equation \eqref{BSprod} is sent to the product $(\alpha\gamma+\beta\delta)+(\alpha\delta+\beta\gamma)\varepsilon$ of our two elements of $\mathbb{E}^{1}$. Therefore our map yields an isomorphism from the projectivization of the semi-group from \cite{[BS]} onto the multiplicative group $\mathbb{E}^{1}$.

As for the structure of $\mathbb{E}^{1}$, the isomorphism $\mathbb{E}\cong\mathbb{Q}\times\mathbb{Q}$ from Lemma \ref{Eprop} implies that $\mathbb{E}^{\times}\cong\mathbb{Q}^{\times}\times\mathbb{Q}^{\times}$, and the formula for $N^{\mathbb{E}}_{\mathbb{Q}}$ in that lemma shows that $\mathbb{E}^{\times}$ indeed consists of those elements whose $N^{\mathbb{E}}_{\mathbb{Q}}$-image does not vanish, and that $\mathbb{E}^{1}$ is the isomorph $\big\{\big(r,\frac{1}{r}\big)\big|r\in\mathbb{Q}^{\times}\big\}$ of $\mathbb{Q}^{\times}$. As the previous paragraph shows that this is also the structure of the (semi-)group in which we are interested, this proves the proposition.
\end{proof}

\begin{rmk}
Few observations arise from the proof of Proposition \ref{projgp}. First, the formula for the product in $\mathbb{E}^{1}$ appearing in the proof of Proposition \ref{projgp}, which holds for arbitrary two elements of $\mathbb{E}$, is the extension of the multiplication rule of $\alpha+\beta\sqrt{d}$ and $\gamma+\delta\sqrt{d}$ from the field $\mathbb{K}=\mathbb{Q}(\sqrt{d})$ to the case $d=1$. Second, the proof of Proposition \ref{projgp} determines the group structure of $\operatorname{Spin}(\mathbb{E})$, since the latter group is $\mathbb{E}^{1}$ by Corollary \ref{Spin} (with $d=1$). This is related to the fact that the action of an arbitrary element $(r,s)\in\mathbb{E}^{\times}$ multiplies the idempotent $\frac{1+\varepsilon}{2}$ (resp. $\frac{1-\varepsilon}{2}$) by $\frac{r}{s}$ (resp. $\frac{s}{r}$), because of the division by the conjugate $(s,r)$, and this operation has spinor norm $rs\in\mathbb{Q}^{\times}/(\mathbb{Q}^{\times})^{2}$. Restricting to elements of $\mathbb{E}^{1}$, with $s=\frac{1}{r}$, this spinor norm is indeed trivial (for other elements of $\mathbb{E}^{\times}$ acting with trivial spinor norms, see Remark \ref{normsq}  below). Finally, the isomorphism $\mathbb{E}^{1}\cong\mathbb{Q}^{\times}$ from the proof of Proposition \ref{projgp} is in line with the structure result for $\mathbb{K}^{1}$ with non-square $d$ above, as seen in Remark \ref{strucK1}: Indeed, the equivalent of the ring of integers in $\mathbb{E}$ is generated as a ring by any one of the idempotents $\frac{1\pm\varepsilon}{2}$, it is isomorphic to $\mathbb{Z}\times\mathbb{Z}$ under the map from Lemma \ref{Eprop}, and every prime ideal of $\mathbb{Z}$ splits there (since there are two distinct ideals in that ``ring of integers'' lying over it). The units in that ``ring of integers'' are just $\pm1$ and $\pm\varepsilon$ (even though $d=1>0$, and for real quadratic fields the group of units has rank 1), producing only the $\pm1$ factor in $\mathbb{Q}^{\times}$ since the units $\pm\varepsilon$ have norm $-1$ and hence do not belong to $\mathbb{E}^{1}$. As one ``ideal'' lying over $p$ is generated by $\pi=(p,1)$, and the quotient $\frac{\pi}{\pi^{\iota}}$ is $\big(p,\frac{1}{p}\big)$, the positive part $\mathbb{Q}^{\times}_{+}$ of $\mathbb{Q}^{\times}$ is indeed generated as in Remark \ref{strucK1}. \label{actofE1}
\end{rmk}

Note that Proposition \ref{projgp} is essentially Theorem 4 of \cite{[M]}, but with the difference that in that reference the projective relation involves only positive scalars (hence the extra multiplier of $\{\pm1\}$ appearing in that reference), and the excluded Pythagorean triples with $a=0$ do appear there, with vanishing $\mathbb{Q}$-images. We also remark that the underlying spinor kernel $\operatorname{SO}^{1}(\mathbb{E})$ from Corollary \ref{Spin} is obtained by dividing $\operatorname{Spin}(\mathbb{E})=\mathbb{E}^{1}$ by $\pm1$ (see, e.g., Equation \eqref{shexseq} with $\mathbb{K}=\mathbb{E}$). But as the latter group is the product of $\{\pm1\}$ and $\mathbb{Q}^{\times}_{+}$ (as in Remark \ref{actofE1}), we deduce that $\operatorname{SO}^{1}(\mathbb{E})\cong\mathbb{Q}^{\times}_{+}$, and the short exact sequence from Equation \eqref{shexseq} splits in this case. This is not always the case, as one sees by considering $\mathbb{K}=\mathbb{Q}(\sqrt{-1})$ and observing that $-1$ is the square of a non-trivial element there, but it does always happen when $\mathbb{K}=\mathbb{Q}(\sqrt{d})$ and $d>1$ (i.e., when $\mathbb{K}$ is a real quadratic field), by working with a real embedding of $\mathbb{K}$.

\smallskip

Allowing the entries of Pythagorean triples to be rational, we deduce the following consequence.
\begin{cor}
The set of all \emph{rational} Pythagorean triples $(a,b,c)$ in which $a\neq0$ forms, with the product from Equation \eqref{BSprod}, a group that is isomorphic to $\mathbb{Q}^{\times}\times\mathbb{Q}^{\times}$. The subset of classical, integral Pythagorean triples $(a,b,c)$ with $a\neq0$, appearing in \cite{[BS]}, is a sub-semi-group of this group. \label{totgp}
\end{cor}

\begin{proof}
Consider the map sending a rational Pythagorean triple $(a,b,c)$ with $a\neq0$ to the pair consisting of $a\in\mathbb{Q}^{\times}$ and the element $\alpha+\beta\varepsilon\in\mathbb{E}^{1}$ appearing in the proof of Proposition \ref{projgp}. This map is bijective, since given $a$ and $\alpha+\beta\varepsilon$ we can reproduce $b$ and $c$ as $a\beta$ and $a\alpha$ respectively. Moreover, the multiplicativity of the $a$-entries in Equation \eqref{BSprod} and the proof of Proposition \ref{projgp} show that when the rational Pythagorean triples are endowed with that product, this map becomes a homomorphism (hence an isomorphism) of semi-groups. Since $\mathbb{Q}^{\times}\times\mathbb{E}^{1}$ is a group, which is isomorphic to $\mathbb{Q}^{\times}\times\mathbb{Q}^{\times}$ by the proof of Proposition \ref{projgp}, this proves the first assertion. As the product in Equation \eqref{BSprod} preserves integrality, the second assertion follows as well. This proves the corollary.
\end{proof}
The projectivization in Proposition \ref{projgp} hints that it might be more natural to consider the group from Corollary \ref{totgp}, not as an isomorph of $\mathbb{Q}^{\times}\times\mathbb{Q}^{\times}$, but rather as lying in the short exact sequence
\begin{equation}
1\to\{(a,0,a)|a\in\mathbb{Q}^{\times}\}\to\{\mathrm{rational\ Pythagorean\ triples\ with\ }a\neq0\}\to\mathbb{E}^{1}\to1 \label{seqtosplit}
\end{equation}
via the projection onto $\mathbb{E}^{1}$. Corollary \ref{totgp} implies that this short exact sequence splits, and indeed, the kernel of the projection to $a$ yields such a splitting using the subgroup $\{(1,\beta,\alpha)|\alpha+\beta\varepsilon\in\mathbb{E}^{1}\}$. However, we shall use a different set-theoretic splitting for the sequence from Equation \eqref{seqtosplit} when we investigate the sub-semi-group of integral Pythagorean triples in Section \ref{Explicit} below.

\begin{rmk}
The product rule from Equation \eqref{BSprod} implies that a (rational) Pythagorean triple $(a,b,c)$ is related to the element $c+b\varepsilon$ of $\mathbb{E}$. Indeed, this element has norm $a^{2}$, and since we require that $a\neq0$, such elements are all invertible. They comprise the subgroup of $\mathbb{E}^{\times}$ defined by the condition that the norm is a square in $\mathbb{Q}^{\times}$. A similar subgroup of $\mathbb{K}^{\times}$ can be considered also for non-square $d$, and in the language of Proposition \ref{GSpin}, Corollary \ref{Spin}, and the short exact sequences, this group is the inverse image of $\operatorname{SO}^{1}(V)$ inside $\operatorname{GSpin}(V)$. The difference, for $d=1$, between this group and the group from Corollary \ref{totgp}, is that while $c+b\varepsilon$ determines $a^{2}$ as its norm, the sign of $a$ is an extra piece of information from a Pythagorean triple $(a,b,c)$. It therefore follows from Corollary \ref{totgp} that this subgroup of $\mathbb{E}^{\times}$ is isomorphic to $\mathbb{Q}^{\times}_{+}\times\mathbb{Q}^{\times}$ (by canonically taking $a$ to be the positive square root of the norm of the element $c+b\varepsilon$). We also note that $|\alpha|>|\beta|$ for any element $\alpha+\beta\varepsilon$ of $\mathbb{E}^{1}$ (and, in fact, for any element of $\mathbb{E}^{\times}$ whose norm is positive), meaning that in the product of this element with another such element $\gamma+\delta\varepsilon$, the sign of the first coefficient $\alpha\gamma+\beta\delta$ is the product of the signs of $\alpha$ and $\gamma$. Moreover, the inverse of $\alpha+\beta\varepsilon$ is the conjugate $\alpha-\beta\varepsilon$ divided by the norm, and the first coordinate has the same sign as $\alpha$ when the norm is positive (in particular for elements of $\mathbb{E}^{1}$). Hence $\mathbb{E}^{1}$ factors as $\{\pm1\}$ times the subgroup $\mathbb{E}^{1}_{+}$ consisting of those elements $\alpha+\beta\varepsilon\in\mathbb{E}^{1}$ in which $\alpha>0$, and this is the subgroup splitting the short exact sequence from Equation \eqref{shexseq} in this case. Such elements come from Pythagorean triples $(a,b,c)$ in which $a$ and $c$ have the same sign---see also Remark \ref{invsign} below. \label{normsq}
\end{rmk}

As mentioned in the Introduction, references like \cite{[E]} and \cite{[T]} put the $c$-coordinate of the Pythagorean triple $(a,b,c)$ in the denominator, thus obtaining the group $\mathbb{K}^{1}$, for the Gaussian field $\mathbb{K}=\mathbb{Q}(\sqrt{-1})$, as a subgroup of the circle group $\big\{z\in\mathbb{C}\big||z|=1\big\}$. As Remark \ref{strucK1} shows, this group is isomorphic, canonically up to local inversions, to the subgroup of $\mathbb{Q}^{\times}_{+}$ generated by the primes $p$ that are congruent to 1 modulo 4 times the torsion group of order 4 generated by $\sqrt{-1}$ (the reference \cite{[E]} essentially divides by the latter torsion subgroup to get a free Abelian group). Some of the results of this section and the following ones have analogues in this point of view as well, but we shall concentrate on the group structure coming from $\mathbb{E}$, which is much more canonical. Indeed, the only choice involved here is the distinction between $\varepsilon$ and $-\varepsilon$, inversion of which means global multiplicative inversion of $r$, and the choice of one copy of $\mathbb{Z}\subseteq\mathbb{Q}$ in Lemma \ref{Eprop} determines the ``prime ideal lying over $p$'' for every prime $p$ simultaneously, while for $\mathbb{K}^{1}$ one must choose a generator for every prime $p$ in $1+4\mathbb{Z}$ independently.

\section{Explicit Expressions \label{Explicit}}

Following \cite{[BS]} and others, we shall be interested, for an element $\alpha+\beta\varepsilon\in\mathbb{E}^{1}$, in the primitive integral Pythagorean triple yielding it. There are two such triples, which differ by a global sign, but taking the $c$-coordinate to be positive yields a canonical choice. The explicit formula for the normalized such Pythagorean triple that is associated with a given element of $\mathbb{Q}^{\times}$ (or of $\mathbb{E}^{1}$), which essentially appears above Corollary 1 of \cite{[M]}, is as follows.
\begin{prop}
Let $r=\frac{m}{n}$ be a non-zero rational number, presented as a reduced quotient of integers. Then the unique primitive integral Pythagorean triple with positive $c$ that is associated with $r$ is $(2mn,m^{2}-n^{2},m^{2}+n^{2})$ in case $mn$ is even and $\big(mn,\frac{m^{2}-n^{2}}{2},\frac{m^{2}+n^{2}}{2}\big)$ if $mn$ is odd. \label{primtrip}
\end{prop}

\begin{proof}
Lemma \ref{Eprop} associates to a given element $\alpha+\beta\varepsilon$ of $\mathbb{E}$ the pair $(r,s)\in\mathbb{Q}\times\mathbb{Q}$ for which this element is $r\frac{1+\varepsilon}{2}+s\frac{1-\varepsilon}{2}$. This yields the formulae $\alpha=\frac{r+s}{2}$ and $\beta=\frac{r-s}{2}$, and equivalently $r=\alpha+\beta$ and $s=\alpha-\beta$. If our element lies in $\mathbb{E}^{1}$, so that $s=\frac{1}{r}$, then we obtain $\alpha=\frac{r}{2}+\frac{1}{2r}$ and $\beta=\frac{r}{2}-\frac{1}{2r}$. Substituting $r=\frac{m}{n}$ (in reduced terms), our element of $\mathbb{E}^{1}$ is $\frac{m^{2}+n^{2}+(m^{2}-n^{2})\varepsilon}{2mn}$. When this element $\alpha+\beta\varepsilon\in\mathbb{E}^{1}$ arises from a Pythagorean triple $(a,b,c)$ as in the proof of Proposition \ref{projgp}, we recall that $\alpha=\frac{c}{a}$ and $\beta=\frac{b}{a}$, so that our Pythagorean triple is $(a,a\beta,a\alpha)$. We need to determine, with our $\alpha$ and $\beta$, the value of $a$ such that this Pythagorean triple is primitive and has positive last entry. Observing that the greatest common divisor of $2mn$, $m^{2}-n^{2}$, and $m^{2}+n^{2}$ is 2 for odd $m$ and $n$ and 1 otherwise, and that $m^{2}+n^{2}>0$ for every $m$ and $n$ that are not both 0, this proves the proposition.
\end{proof}

Note that we do not impose any sign restriction on $m$ and $n$, so that $\frac{-m}{-n}$ yields a reduced presentation of the same number $r$ from Proposition \ref{primtrip}. The formula resulting from this presentation is, however, the same.

\begin{rmk}
Recalling that in a primitive Pythagorean triple $(a,b,c)$, precisely one of $a$ and $b$ is even, we see in Proposition \ref{primtrip} that the parity of $a$ is the parity of $mn$. We also observe that the action of the ``Galois'' automorphism $\iota$ of $\mathbb{E}$ over $\mathbb{Q}$ on an element $\big(r,\frac{1}{r}\big)$ of $\mathbb{E}^{1}$ amounts to inverting $r$ (multiplicatively), hence leaving $\alpha$ invariant and changing the sign of $\beta$. The resulting operation on the Pythagorean triple is the inversion of the sign of $b$, in correspondence with Proposition \ref{primtrip} and the fact that sending $r$ to $\frac{1}{r}$ simply interchanges the integers $m$ and $n$ there. Multiplying $r$ by $-1$ inverts the sign of $a$ alone, and therefore taking $r$ to $-\frac{1}{r}$ corresponds to inverting the signs of both $a$ and $b$. We also deduce from the proof of Proposition \ref{primtrip} that $r>0$ if and only if $\alpha>0$, or equivalently if $a$ and $c$ have the same sign, so that the subgroup $\mathbb{Q}^{\times}_{+}$ of $\mathbb{Q}^{\times}$ indeed corresponds to the subgroup $\mathbb{E}^{1}_{+}$ of $\mathbb{E}^{1}$ appearing in Remark \ref{normsq}. As for the sign of $b$, for $r>0$ it coincides with those of $a$ and $c$ if and only if $r>1$, while if $r<0$ then the sign of $b$ is like that of $a$ when $-1<r<0$, and like that of $c$ if $r<-1$ (this determines the behavior of the involution in the coordinates from Remark \ref{ratinv} below). Finally, there are two primitive Pythagorean triples with positive $c$ that are omitted in Proposition \ref{primtrip}, since their $a$-entries vanish, namely $(0,\pm1,1)$. These are associated, by the same formula, with $r=0=\frac{0}{1}\in\mathbb{Q}$ for the $-$ sign and with $r=\infty=\frac{1}{0}\in\mathbb{P}^{1}(\mathbb{Q})$ when the sign is $+$. The two corresponding elements of an extended version of $\mathbb{E}^{1}$ can both be written as $\infty+\infty\varepsilon$, but for $r=\infty$ these two instances of $\infty$ should be understood to ``have the same sign'', while for $r=0$ they should be viewed as ``additive inverses''. \label{invsign}
\end{rmk}

Considering the usual generating set of $\mathbb{Q}^{\times}$, we immediately deduce the following result.
\begin{cor}
The group from Proposition \ref{projgp} is generated by the elements associated with the Pythagorean triples $(-1,0,1)$ (or equivalently $(1,0,-1)$), $(4,3,5)$, and $\big(p,\frac{p^{2}-1}{2},\frac{p^{2}+1}{2}\big)$ for odd prime $p$. The positive subgroup, associated with $\mathbb{E}^{1}_{+}$ or with $\mathbb{Q}^{\times}_{+}$ via Proposition \ref{projgp} and Remarks \ref{normsq} and \ref{invsign}, is generated by the same set of Pythagorean triples, with the first one excluded. \label{gentrip}
\end{cor}

\begin{proof}
We know that $\mathbb{Q}^{\times}_{+}$ is generated by the prime numbers, and for $\mathbb{Q}^{\times}$ we have to add $-1$ as a generator as well. As when $r$ is $-1$, 2, or an odd prime $p$, Proposition \ref{primtrip} associates the asserted primitive Pythagorean triples, this proves the corollary.
\end{proof}
Note that the prime 2 is separated from the odd primes in Corollary \ref{gentrip} because for $r=2=\frac{2}{1}$ we have to use the formula with even $mn$, while when $p$ is odd, the formula for $r=p=\frac{p}{1}$ that we must apply has odd $mn$. In addition, for elements of $\mathbb{Q}^{\times}_{+}$ the normalizations of taking the primitive Pythagorean triple with $a>0$ (as Corollary \ref{totgp} implicitly suggests) coincides with that in which $c>0$ as in Proposition \ref{primtrip} (since these numbers have the same sign by Remarks \ref{normsq} and \ref{invsign}), while for the other elements we can choose any of these two normalizations, and they are opposites. This is the reason for the two different choices of the representative for the sign generator $-1$ in Corollary \ref{gentrip}.

Because of this issue with the signs, in our analysis of the splitting of the short exact sequence from Equation \eqref{seqtosplit} in a manner that is associated with integral Pythagorean triples and integrality, we shall restrict attention to elements of the subgroup $\mathbb{Q}^{\times}_{+}\times\mathbb{Q}^{\times}_{+}$ in Corollary \ref{totgp}, namely to Pythagorean triples in which both $a$ and $c$ are positive. Hence we take only positive $a$ in the kernel in Equation \eqref{seqtosplit}, and the quotient there reduces to $\mathbb{E}^{1}_{+}$. The full group is obtained as the direct product with the Klein 4-group consisting of the Pythagorean triples $(\pm1,0,\pm1)$ with unrelated signs (the sign of $b$, or more generally of $\frac{b}{a}$, is determined by whether the parameter $r$ is larger or smaller than 1, as Remark \ref{invsign} shows). For doing so we recall that the group from Proposition \ref{projgp} and Corollary \ref{gentrip} is a quotient group of the one from Corollary \ref{totgp}, so that we shall have to say, for every element of the former group, which pre-image (or ``lift'') of it we consider in the latter group. In the direct product from Corollary \ref{totgp} we have taken the ``lift'' with $a=1$, but when we are interested in integral Pythagorean triples, this is not a very good choice. With the restriction $a>0$ and $c>0$, the only natural candidate for that ``lift'' is the primitive integral Pythagorean triple from Proposition \ref{primtrip}. Identifying $\mathbb{E}^{1}_{+}$ with $\mathbb{Q}^{\times}_{+}$ as in the proof of Proposition \ref{projgp} once again, we obtain the following result.
\begin{thm}
Given two pairs $\big(h,\frac{m}{n}\big)$ and $\big(g,\frac{k}{l}\big)$ in $\mathbb{Q}^{\times}_{+}\times\mathbb{Q}^{\times}_{+}$ (with the two fractions being reduced), set their product to be $\big(gh\gcd\{km,ln\}^{2},\frac{km}{ln}\big)$, unless both $mn$ and $kl$ are even, where the product is $\big(4gh\gcd\{km,ln\}^{2},\frac{km}{ln}\big)$ if the reduced form of $\frac{km}{ln}$ involves only odd numbers and $\big(2gh\gcd\{km,ln\}^{2},\frac{km}{ln}\big)$ otherwise. Then the map sending the pair $\big(h,\frac{m}{n}\big)$ to $h$ times the primitive Pythagorean triple from Proposition \ref{primtrip} gives an isomorphism of groups between $\mathbb{Q}^{\times}_{+}\times\mathbb{Q}^{\times}_{+}$ and the subgroup of the group of rational Pythagorean triples from Corollary \ref{totgp} consisting of those Pythagorean triples in which $a$ and $c$ are positive. The sub-semi-group of integral Pythagorean triples with positive $a$ and $c$ is the image of the (multiplicative) sub-semi-group $\mathbb{N}\times\mathbb{Q}^{\times}_{+}$ under this isomorphism. \label{semigrp}
\end{thm}
We remark again that for allowing the signs of $a$ and $c$ to be arbitrary, and obtain the full group $\mathbb{Q}^{\times}_{+}\times\mathbb{Q}^{\times}_{+}$, one simply takes the direct product with the Klein 4-group given above. It is also clear that in our notation the set $\mathbb{N}$ of the natural numbers begins with 1 and not with 0.
\begin{proof}
It is clear (e.g., from the short exact sequence from Equation \eqref{seqtosplit}) that the map in question is bijective. Since our set of Pythagorean triples is a group (by Corollary \ref{totgp} and Remarks \ref{normsq} and \ref{invsign}), once we show that this map respects the operations, it will follow that the given structure on $\mathbb{Q}^{\times}_{+}\times\mathbb{Q}^{\times}_{+}$ is a group structure (which is not obvious from the definition) and that the map is a group isomorphism. It is also clear that this map serves as a splitting of the short exact sequence from Equation \eqref{seqtosplit}, and that the behavior with respect to the first entry is the desired one. The proof of the first assertion thus amounts to taking the primitive Pythagorean triples that are associated with the elements $\frac{m}{n}$ and $\frac{k}{l}$ via Proposition \ref{primtrip}, evaluating their product via Equation \eqref{BSprod}, and showing that it gives the asserted scalar times a primitive Pythagorean triple (which will then be the desired one). Moreover, the fact that our map is a splitting of the short exact sequence from Equation \eqref{seqtosplit} implies that the determination of the scalar can be carried out by checking the $a$-coordinate alone. The result about our sub-semi-group of integral Pythagorean triples, and with it the semi-group of all integral Pythagorean triples from \cite{[BS]}, will then also follow.

We therefore take our fractions $\frac{m}{n}$ and $\frac{k}{l}$, and note that the reduced form of the product $\frac{km}{ln}$ is obtained by canceling $\gcd\{km,ln\}$. Hence we can write the $a$-coordinates of the associated Pythagorean triples from Proposition \ref{primtrip} as $2^{\mu}mn$, $2^{\nu}kl$, and $2^{\kappa}\frac{klmn}{\gcd\{km,ln\}^{2}}$, where $\mu$, $\nu$, and $\kappa$ are elements of $\{0,1\}$ that are determined by the corresponding parity condition in that proposition, and therefore the scalar in question is \[\frac{2^{\mu}mn\cdot2^{\nu}kl}{2^{\kappa}\frac{klmn}{\gcd\{km,ln\}^{2}}}=2^{\mu+\nu-\kappa}\gcd\{km,ln\}^{2}.\] We now observe that if $mn$ and $kl$ are both odd then so is $klmn$ and hence also its reduced version, while if precisely one of $mn$ and $kl$ is even then so is the reduced version of $klmn$. Therefore if not both $mn$ and $kl$ are even (i.e., when not both $\mu$ and $\nu$ are 1), then $\kappa=\mu+\nu$ and the scalar in question is the desired one. On the other hand, if both $mn$ and $kl$ are even then $\mu+\nu=2$, and our distinction is precisely between the two possible values of $\kappa$, yielding the required scalar in both of these cases as well. This completes the proof of the theorem.
\end{proof}
As $\gcd\{km,ln\}$ is the product of $\gcd\{k,n\}$ and $\gcd\{m,l\}$, one can verify in the proof of Theorem \ref{semigrp} directly, also in the coordinates $b$ and $c$, that the product of the Pythagorean triples associated with $\frac{m}{n}$ and with $\frac{k}{l}$ via Proposition \ref{primtrip} is the desired scalar times the Pythagorean triple corresponding to $\frac{km}{ln}$. As an example, we consider the primitive Pythagorean triples associated with $r=\frac{m}{n}$ and $\frac{1}{r}=\frac{n}{m}$, whose product in the quotient group from Proposition \ref{projgp} is trivial. One easily verifies that their product via Equation \eqref{BSprod} (namely in Corollary \ref{totgp}) is not the trivial element $(1,0,1)$, but rather its multiple by the scalar $m^{2}n^{2}$ for odd $mn$ and $4m^{2}n^{2}$ for even $mn$.

\section{The Involution \label{Inv}}

It is clear that if $(a,b,c)$ is a Pythagorean triple then so is $(b,a,c)$, and that the operation taking the former to the latter is an involution. On the other hand, our choice of map to $\mathbb{E}^{1}$, hence to $\mathbb{Q}^{\times}$, distinguishes a triple from its image under the involution. Moreover, we shall have to restrict attention to triples in which $b\neq0$ as well, in order for the $a$-coordinate of the image under the involution not to vanish. The resulting subsets of $\mathbb{E}^{1}$ and $\mathbb{Q}^{\times}$ on which the involution is defined, and the formulae for the involution there, are as follows.
\begin{prop}
The involution is defined on an element $\alpha+\beta\varepsilon\in\mathbb{E}^{1}$ if and only if $\beta\neq0$, and then it takes this element to $\frac{\alpha+\varepsilon}{\beta}$. Similarly, the involution is defined on the set $\mathbb{Q}^{\times}\setminus\{\pm1\}$, and it takes $r$ in that set it its Cayley transform $\frac{r+1}{r-1}$. \label{inv}
\end{prop}

\begin{proof}
The condition $b\neq0$ that is required for the image of $(a,b,c)$ under the involution to have non-vanishing $a$-entry is equivalent (since $a\neq0$) to the non-vanishing of $\beta=\frac{b}{a}$. As interchanging $a$ and $b$ takes this value of $\beta$ to $\frac{1}{\beta}$ and the number $\alpha=\frac{c}{a}$ to $\frac{c}{b}=\frac{\alpha}{\beta}$, this proves the first assertion. We recall from the proof of Proposition \ref{primtrip} that the initial parameter $r$ was $\alpha+\beta$, and note that the excluded values are those in which $\beta=0$, so that $\alpha=\pm1$ for $\alpha+0\varepsilon$ to be in $\mathbb{E}^{1}$. The number that we now have to take is therefore $\frac{\alpha+1}{\beta}$, and substituting the values $\alpha=\frac{r}{2}+\frac{1}{2r}$ and $\beta=\frac{r}{2}-\frac{1}{2r}$, the latter quotient becomes $\frac{r+2+1/r}{r-1/r}$, which yields the desired value after canceling $\frac{r+1}{r}$ (this is possible since we already know that we do this operation only when $r\neq1$). This proves the proposition.
\end{proof}
We remark that the excluded values $\pm1$ in Proposition \ref{inv} are precisely those which are taken by the Cayley transform there to 0 or to $\infty$, which are not in $\mathbb{Q}^{\times}$. However, the extensions of the formulae for Pythagorean triples with vanishing $a$-coordinate and the remaining values in $\mathbb{P}^{1}(\mathbb{Q})$ (as well as some infinite elements extending $\mathbb{E}^{1}$), given at the end of Remark \ref{invsign}, yield valid formulae here as well.

\begin{rmk}
In Proposition \ref{primtrip} we have separated the formula for the primitive Pythagorean triple (with positive $c$-entry) that is associated with the rational number $r=\frac{m}{n}$ according to the parity of the product $mn$. The involution takes a reduced fraction with even product to one with an odd product, and vice versa. Indeed, when $mn$ is even the fraction $\frac{m+n}{m-n}$ is reduced, with odd numerator and denominator. On the other hand, if $mn$ is odd then both $m+n$ and $m-n$ are even and exactly one of them is divisible by 4, so that $\frac{(m+n)/2}{(m-n)/2}$ is reduced and $\frac{m+n}{2}\cdot\frac{m-n}{2}$ is even. This is in correspondence with the formulae from Proposition \ref{primtrip}, since when $mn$ is even we can write the entries $2mn$, $m^{2}-n^{2}$, and $m^{2}+n^{2}$ as $\frac{(m+n)^{2}-(m-n)^{2}}{2}$, $(m+n)(m-n)$, and $\frac{(m+n)^{2}+(m-n)^{2}}{2}$ respectively, while if $mn$ is odd then the entries $mn$, $\frac{m^{2}-n^{2}}{2}$, and $\frac{m^{2}+n^{2}}{2}$ are $\big(\frac{m+n}{2}\big)^{2}-\big(\frac{m-n}{2}\big)^{2}$, $2\frac{m+n}{2}\cdot\frac{m-n}{2}$, and $\big(\frac{m+n}{2}\big)^{2}+\big(\frac{m-n}{2}\big)^{2}$ respectively. \label{parinv}
\end{rmk}

\begin{rmk}
Proposition \ref{inv} expresses the action of the involution in the quotient considered in Proposition \ref{projgp}. For the extension to the full group of rational Pythagorean triples (or more precisely, the subset consisting of those triples in which $b\neq0$), one may ask what happens in the second $\mathbb{Q}^{\times}$-coordinate. This depends on the splitting of the short exact sequence from Equation \eqref{seqtosplit} that we choose to use. With the group-theoretic splitting, which is based on the subgroup of rational Pythagorean triples of the form $(1,\beta,\alpha)$ with $\alpha+\beta\varepsilon\in\mathbb{E}^{1}$, interchanging of coordinates but rescaling to land inside this splitting subgroup takes the latter element (with $\beta\neq0$) to $\big(1,\frac{1}{\beta},\frac{\alpha}{\beta})$, yielding another copy of the formula for the action on $\mathbb{E}^{1}$ appearing in Proposition \ref{inv}. On the other hand, with the set-theoretic splitting from Theorem \ref{semigrp} which is based on primitive integrality, all that can change in this coordinate is the sign (since $(a,b,c)$ is integral and primitive if and only if $(b,a,c)$ is). The actual formula depends on whether the primitive integral lift that we take is the one in which $a>0$ (like in Corollary \ref{totgp}), or the one in which $c>0$ (as in Proposition \ref{primtrip}). It is clear that in the latter convention the second $\mathbb{Q}^{\times}$-coordinate remains always invariant, while in the former one it will change sign in case $\beta<0$ (i.e., when $0<r<1$ or when $r<-1$---see Remark \ref{invsign}). \label{ratinv}
\end{rmk}

Proposition \ref{inv} also yields the following result of \cite{[BS]}.
\begin{cor}
All the Pythagorean triples are generated, up to global scalars, by the triples $(4,3,5)$ and $(-1,0,1)$ (or $(1,0,-1)$) using the multiplication from Equation \eqref{BSprod}, the operation of inverting the sign of $b$, and the involution. \label{geninv}
\end{cor}

\begin{proof}
Proposition \ref{projgp} allows us to investigate the question in the coordinate $r$, and as the sign is covered by multiplication with the last asserted Pythagorean triple, we may restrict attention to Pythagorean triples for which this coordinate is positive. Since the operations in question are the group operations, it suffices to verify that for every prime $p$ we can obtain in this manner a Pythagorean triple whose image in the coordinate $r$ is $p$ (as in Corollary \ref{gentrip}). We work by induction on this set of primes (ordered by magnitude, of course), where the case $p=2$ is immediate because this is the value of $r$ for our initial Pythagorean triple $(4,3,5)$ (see, e.g., Proposition \ref{primtrip}). Assuming that we have generated all the primes that are smaller than the odd prime $p>2$, note that the numerator and the denominator of the Cayley transform $\frac{p+1}{p-1}$ (or more precisely $\frac{(p+1)/2}{(p-1)/2}$, as in Remark \ref{parinv}) of $p$ are divisible only by primes that are smaller than $p$, so that the induction hypothesis implies that we can generate from $(4,3,5)$ a Pythagorean triple with this $r$-value using our operations. The formula for the involution from Proposition \ref{inv} then allows us to generate a Pythagorean triple for which $r=p$ as desired. This proves the proposition.
\end{proof}

We remark that with the multiplication law of \cite{[E]} and \cite{[T]}, which is based on $\mathbb{K}=\mathbb{Q}(\sqrt{-1})$, our involution corresponds to complex conjugation times a power of $\sqrt{-1}$. Dividing by scalars, as well as by the torsion part (as \cite{[E]} does), this involution reduces to the inversion on the free quotient group of $\mathbb{K}^{1}$ from Remark \ref{strucK1}.

\section{Small Differences of Coordinates \label{QuadUnits}}

Recall that \cite{[M]} defines the \emph{height} $h$ of the Pythagorean triple $(a,b,c)$ to be $c-b$. In the primitive case considered in Proposition \ref{primtrip}, associated with $r=\frac{m}{n}$ and with positive $c$, its value is $n^{2}$ when $mn$ is odd, and $2n^{2}$ if it is even. Our parameter $r$ therefore coincides with the quotient $\frac{a}{h}$ (since $a$ equals $mn$ or $2mn$ respectively), yielding the explicit relation between our Proposition \ref{projgp} and Theorem 4 of that reference. Since the height is used in some enumerations of Pythagorean triples (like the \emph{height-excess enumeration} from \cite{[M]} and others), we can equally use the denominator $n$ of our parameter $r$ for such enumerations, at least in the primitive setting. In particular, the (primitive) Pythagorean triples for which this difference is 1 are precisely those of the form $\big(m,\frac{m^{2}-1}{2},\frac{m^{2}+1}{2}\big)$ for an odd integer $m$ (with $r=m$), and the height $h=2$ is obtained precisely for the triples $(2m,m^{2}-1,m^{2}+1)$ with $m$ an even integer (and again $r=m$). For other heights (or denominators) we shall have to distinguish between enumerating Pythagorean triples and enumerating primitive Pythagorean triples, but in any case the numerator $m$, and with it the parameter $r$, will tend to $\infty$ along any such sequence (when one identifies $+\infty$ and $-\infty$, and in the projective line $\mathbb{P}^{1}(\mathbb{Q})$). We remark that the other parameter from that enumeration, namely the excess $e=a+b-c$, equals $n(m-n)$ when $mn$ is odd and $2n(m-n)$ if $mn$ is even (for the Pythagorean triple associated with our $r$ by Proposition \ref{primtrip}). As the maximal number $d$, called the \emph{increment} in \cite{[M]}, whose square divides $2h$, is $n$ or $2n$ respectively, we find that the parameter $k=\frac{e}{d}$ from that reference is just $m-n$.

There are two additional differences of coordinates that one might consider, where for the difference $c-a$ all that we have to do is to apply our involution to the results involving the height. Indeed, considering the presentations from Remark \ref{parinv} for the Pythagorean triple associated with $r=\frac{m}{n}$ by Proposition \ref{primtrip}, one sees that if $mn$ is odd then the difference $c-a$ is twice the square of the denominator $\frac{m-n}{2}$ of the reduced form of the image of $r$ under the involution, while for even $mn$ it is the square the denominator $m-n$ of that reduced image. The Pythagorean triples for which this difference $c-a$ is 1 or 2, which are the images of those from the previous paragraph under the involution, are associated with the value $\frac{(m+1)/2}{(m-1)/2}=1+\frac{1}{(m-1)/2}$ for odd $m$ with difference 1, and with $\frac{m+1}{m-1}=1+\frac{2}{m-1}$ when $m$ is even for the difference 2. The limit value of $r$ here, as well as for any sequence in which $r$ is $\frac{(m+n)/2}{(m-n)/2}$ or $\frac{m+n}{m-n}$, $n$ is fixed, and $m$ goes to $\infty$, is 1, which produces the Pythagorean triple $(1,0,1)$ for which the difference in question vanishes. Similarly, Remark \ref{invsign} attaches the Pythagorean triple $(0,1,1)$ to the value $r=\infty$, so that once again the limit value of $r$ for the sequences with fixed $c-b$ corresponds to the Pythagorean triple in which the difference in question vanishes.

\smallskip

As for the remaining difference $b-a$, it is clear from the definition that it may vanish if and only if the Pythagorean triple in question is a fixed point of our involution. In this case we would have $\beta=1$, and $r$ is a solution to any of the two equations $\frac{r}{2}-\frac{1}{2r}=1$ or $r=\frac{r+1}{r-1}$, both of which are equivalent to the quadratic equation $r^{2}-2r-1=0$. The real solutions of that equation are $r=1\pm\sqrt{2}$, which are no longer rational, but by extending scalars in $\mathbb{E}$ this parameter indeed corresponds to the norm 1 element $\pm\sqrt{2}+\varepsilon$ as well as to the ``real Pythagorean triple'' $(1,1,\pm\sqrt{2})$. Hence one may expect that integral Pythagorean triples in which this difference is small would have $r$-values that approximate the roots of this equation.

This indicates that the difference $b-a$, as well as some formulae associated with the involution, are related to the real quadratic field $\mathbb{K}=\mathbb{Q}(\sqrt{2})$. Indeed, we recall from Remark \ref{parinv} that if $\frac{m}{n}$ and $\frac{k}{l}$ are two reduced fractions that are related by the involution then the parities of $mn$ and $kl$ are different, and we can prove the following result.
\begin{lem}
Let $\frac{m}{n}$ be $\frac{k}{l}$ reduced fractions with $kl$ odd. Then these fractions are related by our involution if and only if the equality
\begin{equation}
k+m\sqrt{2}=(1+\sqrt{2})(l+n\sqrt{2}) \label{prodinQsq2}
\end{equation}
holds. When this is the case, write $(a,b,c)$ and $(b,a,c)$ for the associated primitive integral Pythagorean triples from Proposition \ref{primtrip}, and then the difference $b-a$ coincide with the norm $N^{\mathbb{K}}_{\mathbb{Q}}(l+n\sqrt{2})$. \label{unitinv}
\end{lem}

\begin{proof}
Since $mn$ is even, we recall from Remark \ref{parinv} that if the involution takes $\frac{m}{n}$ to $\frac{k}{l}$ then $l=m-n$ (so that $m=n+l$) and $k=m+n=2n+l$, which implies Equation \eqref{prodinQsq2}. Conversely, if $mn$ is even then we can write $m=n+l$ for odd $l$, and as Equation \eqref{prodinQsq2} yields also $k=2n+l$, we deduce that $kl$ is odd and the values of $k$ and $l$ in terms of $m$ and $n$ show that $\frac{m}{n}$ and $\frac{k}{l}$ are taken to one another by the involution. Assuming that this is the case, we apply Proposition \ref{primtrip} to the fraction $\frac{m}{n}=\frac{n+l}{n}$ (with even $mn$), and obtain that the difference between $b=(n+l)^{2}-n^{2}=2ln+l^{2}$ and $a=2n(l+n)$ is indeed the asserted norm $l^{2}-2n^{2}$. This proves the lemma.
\end{proof}

\begin{rmk}
Since the multiplier $1+\sqrt{2}$ from Equation \eqref{prodinQsq2}, which is also the positive solution for the equation for $r$ arising from $a=b$, has norm $-1$, we deduce from Lemma \ref{unitinv} that the norm $N^{\mathbb{K}}_{\mathbb{Q}}(k+m\sqrt{2})$ of the other expression from that equation coincides with $a-b$. Moreover, the fraction $\frac{k}{l}=\frac{2n+l}{l}$ indeed yields via Proposition \ref{primtrip} the entries $l(2n+l)=b$ and $\frac{(2n+l)^{2}-l^{2}}{2}=a$, and the same $c$-entry $\frac{(2n+l)^{2}+l^{2}}{2}=2n^{2}+2nl+l^{2}=(n+l)^{2}+n^{2}$. Note that we assume nothing about the signs of $k$, $l$, $m$, and $n$ here, and that the cases in which $n=0$ or $m=0$ are also allowed (with the extension from Remark \ref{invsign}), where the two sides of Equation \eqref{prodinQsq2} become $1+\sqrt{2}$ and $1=(1+\sqrt{2})(-1+\sqrt{2})$ respectively.
\label{altdesc}
\end{rmk}
It follows from Lemma \ref{unitinv} that finding primitive Pythagorean triples $(a,b,c)$ in which $b-a\in\{\pm1\}$ is equivalent to solving Pell's equation for our field $\mathbb{K}=\mathbb{Q}(\sqrt{2})$. These Pythagorean triples are thus characterized as follows.
\begin{thm}
Given $j\in\mathbb{Z}$, write the number $(1+\sqrt{2})^{j}$ as $s_{j}+t_{j}\sqrt{2}$ with integral $s_{j}$ and $t_{j}$. Then the Pythagorean triples $(a,b,c)$ that satisfy $c>0$ and $b-a\in\{\pm1\}$ are precisely those that are associated with the quotients $\frac{t_{j+1}}{t_{j}}$ and with the quotients $\frac{s_{j+1}}{s_{j}}$, and these numbers for the same $j$ are related by the involution. Moreover, the limit of both quotients as $j\to\infty$ is $1+\sqrt{2}$. \label{abclose}
\end{thm}

\begin{proof}
The number $1+\sqrt{2}$ from Equation \eqref{prodinQsq2} is the fundamental unit of the ring of integers in $\mathbb{K}$ (the other solution for $r$ is minus its inverse) and its norm is $-1$. Since our field is real quadratic, any integral element of it that has norm $\pm1$ is of the form $\pm(1+\sqrt{2})^{j}$ for some $j\in\mathbb{Z}$. As the resulting number is a unit, it is not divisible by $\sqrt{2}$ and hence $s_{j}$ is odd for every $j$ (this can also be proved by induction since $s_{0}=1$ and $s_{j+1}=s_{j}+2t_{j}$ for every $j$). We therefore substitute this number to be $l+n\sqrt{2}$ in Equation \eqref{prodinQsq2}, and we may ignore the external sign since it will not affect any of the resulting fractions. Now, with $l=s_{j}$ and $n=t_{j}$ this equation yields $k=s_{j+1}$ and $m=t_{j+1}$, and the fact that $t_{j+1}=s_{j}+t_{j}$ with $s_{j}$ odd shows that indeed $mn$ is even and $kl$ is odd. The first two assertions thus follow from Lemma \ref{unitinv} (the choice of whether the associated quotient is $\frac{t_{j+1}}{t_{j}}$ or $\frac{s_{j+1}}{s_{j}}$ depends on whether our Pythagorean triple has even $a$ and difference $b-a=(-1)^{j}$ or odd $a$ and $a-b=(-1)^{j}$, but we have seen that the involution takes each one to the other). As for the limit, note that conjugating the definition implies that $s_{j}-t_{j}\sqrt{2}=(1-\sqrt{2})^{j}$, so that $s_{j}$ and $t_{j}$ can be presented as $\frac{(1+\sqrt{2})^{j}+(1-\sqrt{2})^{j}}{2}$ and $\frac{(1+\sqrt{2})^{j}-(1-\sqrt{2})^{j}}{2\sqrt{2}}$ respectively. Since our expression for $s_{j}-t_{j}\sqrt{2}$ decreases exponentially as $j\to\infty$, the limits of the two quotients are as desired. This proves the theorem.
\end{proof}
For the first few examples of Theorem \ref{abclose}, we evaluate the next few powers of $1+\sqrt{2}$ as $3+2\sqrt{2}$, $7+5\sqrt{2}$, $17+12\sqrt{2}$, $41+29\sqrt{2}$, and $99+70\sqrt{2}$. The respective quotients $\frac{t_{j+1}}{t_{j}}$ are 2, $\frac{5}{2}$, $\frac{12}{5}$, $\frac{29}{12}$, and $\frac{70}{29}$, producing the Pythagorean triples $(4,3,5)$, $(20,21,29)$, $(120,119,169)$, $(696,697,985)$, and $(4060,4059,5741)$. The quotients $\frac{s_{j+1}}{s_{j}}$ are 3, $\frac{7}{3}$, $\frac{17}{7}$, $\frac{41}{17}$, and $\frac{99}{41}$, always with odd numbers, and one can verify that they are the respective images of the previous numbers under the involution and that they produce the same Pythagorean triples with $a$ and $b$ inverted. The triple $(0,1,1)$ and its involutive image $(1,0,1)$, associated with $r=\infty$ and $r=1$ respectively via Remark \ref{invsign}, correspond to $j=0$, and since the inverse of $1+\sqrt{2}$ is $-1+\sqrt{2}$, we get $t_{j}=(-1)^{j}t_{-j}$ and $s_{j}=(-1)^{j-1}s_{-j}$ for any $j$. Hence the values arising from $-j-1$ are obtained from those associated with $j$ by the map sending $r$ to $-\frac{1}{r}$, which corresponds to inverting the signs of both $a$ and $b$ (see Remark \ref{invsign}). Therefore the limit as $j\to-\infty$ in Theorem \ref{abclose} is the second solution $-\frac{1}{1+\sqrt{2}}=1-\sqrt{2}$ of the quadratic equation defining the fundamental unit $1+\sqrt{2}$.

\begin{rmk}
Ordering (primitive) Pythagorean triples according to the difference $b-a$ was seen in Lemma \ref{unitinv} to depend on norms from $\mathbb{K}$, and the next odd integer that can be obtained as such a norm is 7. Indeed, the element $3-\sqrt{2}$ has norm 7, and the proofs of Lemma \ref{unitinv} and Theorem \ref{abclose} produce a pair of sequences of primitive Pythagorean triples, related via the involution, where in each such Pythagorean triple we have $b-a\in\{\pm7\}$. Moreover, the $r$-values in both sequences tend to the same limit $1+\sqrt{2}$ as $j\to\infty$. In addition, starting with the element $-1+2\sqrt{2}$, of norm $-7$, yields an independent pair of such sequences, with the same properties. As Remark \ref{invsign} shows that replacing $r$ by $\frac{1}{r}$ or by $-r$ inverts the sign of only one of $a$ and $b$, the same analysis investigates Pythagorean triples with small value of $|b+a|$. Similarly, in the Pythagorean triples associated with $r=\frac{1}{m}$ (tending to 0 as $m\to\infty$) the sum $c+b$ is 1 for odd $m$ and 2 for even $m$, where if $r=-1-\frac{2}{m-1}$ (with limit $-1$) we get Pythagorean triples in which $c+a$ equals 1 or 2 (according to the parity of $m$ again). \label{N7andadd}
\end{rmk}

\section{More General Linear Combinations \label{LinComb}}

The height $h=c-b$ from \cite{[M]}, its involutive analogue $c-a$, and the difference $b-a$ considered in Lemma \ref{unitinv} and Theorem \ref{abclose} are the simplest linear expressions in the entries of the Pythagorean triple $(a,b,c)$, and Remark \ref{N7andadd} considers their additive analogues. In this section we investigate relations between more general linear expressions in these entries and additional Pell's equations.

\smallskip

The fundamental unit $1+\sqrt{2}$ played a double role in the analysis of almost fixed points of the involution: One is as the limit of the $r$-values (namely the solution of the quadratic equation yielding the ``real Pythagorean triple'' of the fixed point, and another one is in Equation \eqref{prodinQsq2} producing the relations between involution images. We will obtain a similar double role of a unit in linear expressions like $c-qa$ and $b-qa$ for $q\in\mathbb{N}$, but not with the (reduced) expressions $pc-qa$ and $pb-qa$ with non-trivial (co-prime) $p$ and $q$.

We have seen in the proof of Lemma \ref{unitinv} that when one investigates the difference $b-a$, the more useful presentation of a fraction $\frac{m}{n}$ with even $mn$ is as $\frac{n+l}{n}$ for $l=m-n$. This generalizes to other relations as follows.
\begin{lem}
Let $q\geq2$ be an integer, denote the number $q^{2}\pm1$ by $D_{q,\pm}$, and set $\mathbb{K}_{q,\pm}=\mathbb{Q}\big(\sqrt{D_{q,\pm}}\big)$ for every sign $\pm$. Take a fraction $\frac{m}{n}$ with even $mn$, write it as $\frac{l+qn}{n}$ for some $l$, and let $(a,b,c)$ be the associated Pythagorean triple from Proposition \ref{primtrip}. Then the differences $b-qa$ and $c-qa$ coincide with the norms $N^{\mathbb{K}_{q,+}}_{\mathbb{Q}}\big(l+n\sqrt{D_{q,+}}\big)$ and $N^{\mathbb{K}_{q,-}}_{\mathbb{Q}}\big(l+n\sqrt{D_{q,-}}\big)$ respectively. \label{normp1}
\end{lem}

\begin{proof}
With $m=l+qn$, the fact that $mn$ is even implies that $a=2n(l+qn)$, $b=l^{2}+2qnl+(q^{2}-1)n^{2}$, and $c=l^{2}+2qnl+(q^{2}+1)n^{2}$. The two equalities now follow from simple algbera. This proves the lemma.
\end{proof}

\begin{rmk}
The norm assertion in Lemma \ref{unitinv} is the case with $q=1$ and the sign $+$ in Lemma \ref{normp1}, since $D_{1,+}=2$. On the other hand, we have $D_{1,-}=0$, and in this case we indeed have $c-a=l^{2}$. This is in correspondence with the fact that Lemma \ref{unitinv} determines the image of $\frac{n+l}{n}$ under the involution as $\frac{n+2l}{l}$, and the height of the associated Pythagorean triple is indeed $l^{2}$. As the relation between $m$, $n$, and $l$ here looks like the one arising from Equation \eqref{prodinQsq2} but with $1+\sqrt{2}$ replaced by $q+\sqrt{D_{q,\pm}}$, Remark \ref{altdesc} implies that one should consider the additional fraction $\frac{k}{l}$, where $k$ is obtained by this modified equation, and we write the associated Pythagorean triple from Proposition \ref{primtrip} as $(f,g,h)$. This value of $k$ is $nD_{q,\pm}+ql$, but we have to make the fraction $\frac{k}{l}$ reduced. In addition to $l$ maybe having a common divisor with $D_{q,\pm}$, we note that for odd $q$ (with even $D_{q,\pm}$) the fact that $mn$ is even implies that $l$ is odd hence so is $k$, but if $q$ is even (and $D_{q,\pm}$ is odd) then $l$ must have the opposite parity from $n$ in order for $\frac{m}{n}$ to be reduced with even $mn$. Checking the possible cases shows that the reduced version of $\frac{k}{l}$ has even product $kl$. It follows that for odd $q$ the appropriate difference (namely $g-qf$ for the $+$ sign and $h-qf$ when the sign is $-$) equals \[\frac{N^{\mathbb{K}_{q,\pm}}_{\mathbb{Q}}\big(nD_{q,\pm}+ql\sqrt{D_{q,\pm}}\big)}{2\gcd\{l,D_{q,\pm}\}^{2}}=
\frac{-D_{q,\pm}}{2\gcd\{l,D_{q,\pm}\}^{2}}N^{\mathbb{K}_{q,\pm}}_{\mathbb{Q}}\big(l+n\sqrt{D_{q,\pm}}\big)\] for odd $q$, and just \[\frac{N^{\mathbb{K}_{q,\pm}}_{\mathbb{Q}}\big(nD_{q,\pm}+ql\sqrt{D_{q,\pm}}\big)}{\gcd\{l,D_{q,\pm}\}^{2}}=
\frac{-D_{q,\pm}}{\gcd\{l,D_{q,\pm}\}^{2}}N^{\mathbb{K}_{q,\pm}}_{\mathbb{Q}}\big(l+n\sqrt{D_{q,\pm}}\big)\] when $q$ is even. We remark that when the initial fraction $\frac{m}{n}$ has odd $mn$, the difference from Lemma \ref{normp1} would be half that norm, and then for $\frac{k}{l}$ the product $kl$ would be odd as well for even $q$ (with the difference associated with $\frac{k}{l}$ having 2 in the denominator as well). For odd $q$ we will have even $l$ in this case, and the parity of the reduced form of $kl$ will depend on the parity of $\frac{D_{q,\pm}}{\gcd\{l,D_{q,\pm}\}}$ (the latter is always odd for the $+$ sign, but may be either even or odd when the sign is $-$). \label{klq}
\end{rmk}

Note that for $q\geq2$ the involution plays no role. Indeed, while when $q=1$ and the sign is $+$ the involution only inverted the sign of $b-a$, here it would take Pythagorean triples for which we evaluate $b-qa$ or $c-qa$ to Pythagorean triples with information about the differences $a-qb$ or $c-qb$. The remaining differences $a-qc$ and $b-qc$ are not interesting for $q\geq2$, since $c$ is larger in absolute value from both $a$ and $b$.

\smallskip

The parameter $q+\sqrt{D_{q,\pm}}$ from Remark \ref{klq} also gives the construction of all the values $r\in\mathbb{Q}^{\times}$ for which the difference $b-qa$ or $c-qa$ is $\pm1$, in the following generalization of Theorem \ref{abclose} for such relations with $q\geq2$.
\begin{thm}
Write, for $j\in\mathbb{Z}$, the number $\big(q+\sqrt{D_{q,\pm}}\big)^{j}$ as $s_{q,j}^{\pm}+t_{q,j}^{\pm}\sqrt{D_{q,\pm}}$, where $s_{q,j}^{\pm}$ and $t_{q,j}^{\pm}$ are integers. Then for every $j$ the quotient $\frac{t_{q,j+1}^{+}}{t_{q,j}^{+}}$ represents a Pythagorean triple $(a,b,c)$ with $c>0$ and such that $b-qa\in\{\pm1\}$, while the Pythagorean triple that is associated with $\frac{t_{q,j+1}^{-}}{t_{q,j}^{-}}$ has $c>0$ and $c-qa=1$. The limit of these numbers as $j\to\infty$ is $q+\sqrt{D_{q,\pm}}$ (with the respective sign). Moreover, when $q\geq2$ this produces all the Pythagorean triples $(a,b,c)$ with $c>0$ and such that either $b-qa$ or $c-qa$ are in $\{\pm1\}$, except for the case with $q=2$ and the $-$ sign, where there exist Pythagorean triples $(a,b,c)$ with $c-2a=-1$. The latter also lie in a sequence, which tends with the parameter $j$ (oriented appropriately) to $q+\sqrt{D_{q,-}}=2+\sqrt{3}$. \label{linrel1}.
\end{thm}

\begin{proof}
First we claim that $t_{j,q}^{\pm}t_{j+1,q}^{\pm}$ is even for every $j$. Indeed, we have
\begin{equation}
t_{j+1,q}^{\pm}=qt_{j,q}^{\pm}+s_{j,q}^{\pm}\qquad\mathrm{and}\qquad s_{j+1,q}^{\pm}=qs_{j,q}^{\pm}+D_{q,\pm}t_{j,q}^{\pm} \label{relsqD}
\end{equation}
for every $j$, starting with $s_{1,q}^{\pm}=q$ and $t_{1,q}^{\pm}=1$ (or with $s_{0,q}^{\pm}=1$ and $t_{0,q}^{\pm}=0$), and we know that if $q$ is even then $D_{q,\pm}$ is odd while for odd $q$ we have even $D_{q,\pm}$. It follows that for odd $q$ all the numbers $s_{j,q}^{\pm}$ are odd and hence $t_{j,q}^{\pm}$ has the same parity as $j$, while if $q$ is even then $t_{j,q}^{\pm}$ once again has the parity of $j$ while $s_{j,q}^{\pm}$ is of the opposite parity. As Equation \eqref{relsqD} shows that when writing $\frac{m}{n}=\frac{t_{q,j+1}^{\pm}}{t_{q,j}^{\pm}}$ as $\frac{l+qn}{n}$ we get $l=s_{q,j}^{\pm}$, Lemma \ref{normp1} shows that the appropriate difference (namely $b-qa$ for the $+$ sign and $c-qa$ when the sign is $-$) is $N^{\mathbb{K}_{q,+}}_{\mathbb{Q}}\big(s_{q,j}^{\pm}+t_{q,j}^{\pm}\sqrt{D_{q,\pm}}\big)$. But the number $q+\sqrt{D_{q,\pm}}$ has norm $\mp1$, so that from the multiplicativity of the norm we deduce that the difference in question is $(\mp1)^{j}$, and the asserted numbers indeed produce Pythagorean triples with the required properties. In addition, an argument similar to the proof of Theorem \ref{abclose} allows us to write \[2s_{j,q}^{\pm}=(q+\sqrt{D_{q,\pm}})^{j}+(q-\sqrt{D_{q,\pm}})^{j}\] and \[2\sqrt{D_{q,\pm}} \cdot t_{j,q}^{\pm}=(q+\sqrt{D_{q,\pm}})^{j}-(q-\sqrt{D_{q,\pm}})^{j},\] with $s_{q,j}^{\pm}-t_{q,j}^{\pm}\sqrt{D_{q,\pm}}=\big(q-\sqrt{D_{q,\pm}}\big)^{j}$ decreasing exponentially, from which the assertion about the limit also follows.

It remains to show that these are the only such Pythagorean triples, which by Lemma \ref{normp1} and the part of Remark \ref{klq} considering odd $mn$ is equivalent to showing that all the elements $l+n\sqrt{D_{q,\pm}}$ whose norm is either $\pm1$ or $\pm2$ are the ones we have just presented, unless $q=2$, the sign is $-$, and $D_{2,-}=3$. Now, elements with norm $\pm1$ are invertible elements of the ring $\mathbb{Z}\big[\sqrt{D_{q,\pm}}\big]$, in which it is easy to verify that our element $q+\sqrt{D_{q,\pm}}$ is a fundamental unit (note that this is true in the ring $\mathbb{Z}\big[\sqrt{D_{q,\pm}}\big]$, and usually not in the full ring of integers of $\mathbb{K}_{q,\pm}$). Therefore (up to multiplying by $-1$) all the elements of norm $\pm1$ are our powers $s_{q,j}^{\pm}+t_{q,j}^{\pm}\sqrt{D_{q,\pm}}$ (and in particular there are no elements of norm $-1$ when the sign is $-$, a fact that is also evident since the residue of $D_{q,-}$ modulo 4 is either 3 or 0). Assume now that $N^{\mathbb{K}_{q,\delta}}_{\mathbb{Q}}\big(l+n\sqrt{D_{q,\delta}}\big)=\pm2$ for some sign $\delta\in\{\pm\}$, and we may assume w.l.o.g. that $n$ and $l$ are both positive. Simple algebra shows that \[\big(l-qn-\tfrac{\delta}{2}n\big)\big(l+qn+\tfrac{\delta}{2}n\big)=l^{2}-\big(qn+\tfrac{\delta}{2}n\big)^{2}=l^{2}-q^{2}n^{2}-\delta n^{2}+\tfrac{n^{2}}{4}=\pm2+\tfrac{n^{2}}{4}\] (because of the norm condition), and if $n\geq3$ then the right hand side is smaller than $\frac{n^{2}}{2}$ in absolute value. We divide by $l+qn+\tfrac{\delta}{2}n$ and note that this multiplier is larger than $qn+\tfrac{\delta}{2}n$, and deduce that the second multiplier on the left hand side is smaller than $\frac{n}{2q+\delta}$ in absolute value and therefore $|l-qn|<n$. As multiplying our element $l+n\sqrt{D_{q,\delta}}$, of norm $\pm2$, by the element $q-\sqrt{D_{q,\delta}}$, of norm $-\delta$, yields $(lq-nD_{q,\delta})-(l-qn)\sqrt{D_{q,\delta}}$, we deduce that the existence of elements with the desired norm and with $n\geq3$ implies the existence of such elements with smaller $n$. Moreover, we consider only odd $n$ (for $n(l+qn)$ to be odd, or simply since $n^{2}D_{q,\delta}\pm2$ has residue 2 modulo 4 if $n$ is even, and hence cannot be a square $l^{2}$), so that $n=1$ and we thus look for $q\geq2$ and $\delta\in\{\pm1\}$ such that $D_{q,\delta}\pm2$ is a square. As the latter number is either $q^{2}\pm1$ or $q^{2}\pm3$, this happens only when $q=2$, $\delta=-$, the norm is $-2$, and $l=n=1$.

This shows that no other Pythagorean triple $(a,b,c)$ has $b-qa\in\{\pm1\}$ if $q\geq2$ or $c-qa\in\{\pm1\}$ when $q\geq3$, and it remains to see which additional Pythagorean triples may have $c-2a\in\{\pm1\}$. But our analysis has found the norm $-2$ element $1+\sqrt{3}$, and showed that every other element with such norm can be obtained, up to sign, either from it or from its conjugate by a power of the fundamental unit $2+\sqrt{3}$ (here $\mathbb{Z}\big[\sqrt{D_{2,-}}\big]=\mathbb{Z}[\sqrt{3}]$ is the full ring of integers in $\mathbb{K}_{2,-}=\mathbb{Q}(\sqrt{3})$). Moreover, as $\frac{1+\sqrt{3}}{-1+\sqrt{3}}=2+\sqrt{3}$ as well, we can restrict attention to only one of them. It follows that if we denote $(-1+\sqrt{3})(2+\sqrt{3})^{j}$ by $\xi_{j}+\eta_{j}\sqrt{3}$ (so that $\xi_{j}$ is $3t_{j,2}^{-}-2s_{j,2}^{-}$ and $\eta_{j}$ is the odd number $s_{j,2}^{-}-t_{j,2}^{-}$, or $\xi_{j}=3t_{j-1,2}^{-}+2s_{j-1,2}^{-}$ and $\eta_{j}=s_{j-1,2}^{-}+t_{j-1,2}^{-}$) then the only remaining Pythagorean triples $(a,b,c)$ for which $c-2a\in\{\pm1\}$ are those associated with $\frac{\eta_{j+1}}{\eta_{j}}$, where $\eta_{j}\eta_{j+1}$ is odd (and then $c-2a=-1$). As a small modification of the argument above shows that this quotient also tends to $2+\sqrt{3}$, this completes the proof of the theorem.
\end{proof}

Note that for $j=0$ we obtain the Pythagorean triple $(0,1,1)$ associated with $r=\frac{1}{0}$ by Remark \ref{invsign} for every $q$ and sign $\pm$, in which indeed both $b-qa$ and $c-qa$ are 1. With $j=1$ we get $r=2q$, with the Pythagorean triple $(4q,4q^{2}-1,4q^{2}+1)$ having the two desired properties. Since the first few norm $-2$ elements of the form $\xi_{j}+\eta_{j}\sqrt{3}$ from the proof of Theorem \ref{linrel1} are $-1+\sqrt{3}$, $1+\sqrt{3}$, and $5+3\sqrt{3}$, we obtain from the corresponding ratios $r=1$ and $r=3$ the Pythagorean triples $(1,0,1)$ and $(3,4,5)$ respectively, indeed with $c-2a=-1$.

\begin{rmk}
As $N^{\mathbb{K}_{q,\pm}}_{\mathbb{Q}}\big(q+\sqrt{D_{q,\pm}}\big)=\mp1$, we deduce that $s_{j,q}^{\pm}=(\mp1)^{j}s_{-j,q}^{\pm}$ and $t_{j,q}^{\pm}=(\pm1)^{j}t_{-j,q}^{\pm}$, so that taking $j$ to $-j-1$ replaces $r$ by $\mp\frac{1}{r}$. This is in correspondence with Remark \ref{invsign}, which states that this operation inverts the signs of both $a$ and $b$ (hence that of $b-qa$) when the sign is $+$, while only the sign of $b$ is inverted (and $c-qa$ remain invariant) if the sign is $-$. Other inversions will give Pythagorean triples in which sums like $c+qa$ and $b+qa$ are small, as in Remark \ref{N7andadd}. Note again that unlike the case $q=1$ from Remark \ref{altdesc}, the fractions $\frac{k}{l}$ from Remark \ref{klq} will never give Pythagorean triples with difference $\pm1$ with $q\geq2$ (since $\gcd\{l,D_{q,+}\}$ will be 1 and $D_{q,\pm}>2$). In addition, the parameter $q+\sqrt{D_{q,\pm}}$ is the $r$-value of the ``real Pythagorean triple'' in which the difference in question vanishes, namely $\big(1,q,\sqrt{D_{q,+}}\big)$ for $b-qa$ and $\big(1,\sqrt{D_{q,-}},q\big)$ for $c-qa$, which are related to the elements $\sqrt{D_{q,+}}+q\varepsilon$ or $q+\sqrt{D_{q,-}}\cdot\varepsilon$ of the extended version of $\mathbb{E}^{1}$. As for the involution, it takes our real $r$ (and our elements of the real $\mathbb{E}^{1}$) to values with other properties, namely $\frac{1+\sqrt{D_{q,+}}}{q}$ (and $\frac{\sqrt{D_{q,+}}+\varepsilon}{q}$) for the $+$ sign and $\frac{\sqrt{D_{q,-}}}{q-1}=\sqrt{\frac{q+1}{q-1}}$ (and $\frac{q+\varepsilon}{\sqrt{D_{q,-}}}$) when the sign is $-$. Applying it to the results of Theorem \ref{linrel1} determines the Pythagorean triples in which differences $a-qb$ and $c-qb$ are in $\pm1$, and by appropriate inversions we obtain Pythagorean triples with small $a+qb$ and $c+qb$. \label{addq}
\end{rmk}

\smallskip

We conclude by indicating the partial results that one can get for the more general differences $pb-qa$ and $pc-qa$ when $p$ and $q$ are co-prime and $p>1$. While Remark \ref{addq} will deal with the case $q=1$ for $pb-qa$, we shall need $q>p$ for the difference $pc-qa$ to be interesting (having the involution at hand, the results for $pb-qa$ with $q>p$ also suffice). We then get the following small modification of Lemma \ref{normp1}.
\begin{lem}
Consider a reduced fraction $\frac{m}{n}$ with even $mn$ and such that $p$ divides the denominator $n$, and write it as $\frac{l+qu}{pu}$ for appropriate integers $l$ and $q$. If $(a,b,c)$ is the associated Pythagorean triple from Proposition \ref{primtrip} then $b-\frac{q}{p}a$ is an integer, which coincides with the norm $N^{\mathbb{K}_{+}}_{\mathbb{Q}}(l+u\sqrt{q^{2}+p^{2}})$, for $\mathbb{K}_{+}=\mathbb{Q}(\sqrt{q^{2}+p^{2}})$. In addition, if $q>p$ then $c-\frac{q}{p}a$ is also integral, and its value is $N^{\mathbb{K}_{-}}_{\mathbb{Q}}(l+u\sqrt{q^{2}-p^{2}})$, where $\mathbb{K}_{-}=\mathbb{Q}(\sqrt{q^{2}-p^{2}})$. \label{normgen}
\end{lem}
The proof is just a calculation like in the proof of Lemma \ref{normp1}, combined with the observation that $p$ divides the multiplier $n=pu$ appearing in the formula for $a$ in Proposition \ref{primtrip}.

The determination of the Pythagorean triples in which $pb-qa$ or $pc-qa$ is small encounters several difficulties for general $p$ and $q$. First, Lemma \ref{normgen} deals only with the case in which $p$ divides the denominator $n$ of the associatd fraction $r=\frac{m}{n}$. Second, it might happen that $(p,q,y)$ or $(p,z,q)$ will be Pythagorean for some integers $y$ or $z$, and then $\mathbb{K}_{+}$ or $\mathbb{K}_{-}$ will not be fields. In these cases $pb-qa$ or $pc-qa$ does indeed vanish for some integral Pythagorean triples, and no non-trivial element of the form $N^{\mathbb{K}_{+}}_{\mathbb{Q}}(l+u\sqrt{q^{2}+p^{2}})=l^{2}-u^{2}y^{2}$ or $N^{\mathbb{K}_{-}}_{\mathbb{Q}}(l+u\sqrt{q^{2}-p^{2}})=l^{2}-u^{2}z^{2}$ will equal 1. But also in the field case, the element $q+\sqrt{q^{2} \pm p^{2}}$ expressing $l+qu$ in terms of the coefficients of $l+u\sqrt{q^{2} \pm p^{2}}$ is no longer invertible in the corresponding ring, but has norm $\mp p^{2}$ (this is in correspondence with $\frac{l+qu}{p}$ arising from the product of an element of norm $\pm1$, but this element is no longer integral). Lemma \ref{normgen} will allow us to determine all the Pythagorean triples $(a,b,c)$ in which $p$ divides $a$ and the differences $pb-qa$ or $pc-qa$ is $\pm p$ (provided that the corresponding algebra $\mathbb{K}_{+}$ or $\mathbb{K}_{-}$ is a field, otherwise no such non-trivial Pythagorean triples exist), but for doing so explicitly one will have to find a fundamental unit, which will have nothing to do with our element $q+\sqrt{q^{2} \pm p^{2}}$. Since results where $pb-qa$ or $pc-qa$ is not divisible by $p$ cannot be handled in this method, and also those who do require further calculations, we content ourselves with the results of Theorem \ref{linrel1}, together with the extension of the special case appearing in Theorem \ref{abclose}.

\noindent\textsc{Einstein Institute of Mathematics, the Hebrew University of Jerusalem, Edmund Safra Campus, Jerusalem 91904, Israel}

\noindent E-mail address: zemels@math.huji.ac.il

\end{document}